\documentclass{article} 
\usepackage{amsmath}
\usepackage{amssymb}
\usepackage{amsthm}
\usepackage{multicol}
\usepackage{wrapfig}
\usepackage[retainorgcmds]{IEEEtrantools}

\usepackage{xy}
\input xy
\xyoption{all}

\DeclareMathOperator{\Q}{\mathbb Q}

\DeclareMathOperator{\F}{\mathbb F}

\DeclareMathOperator{\GL}{GL}

\DeclareMathOperator{\PSL}{PSL}

\DeclareMathOperator{\md}{MaxDim}

\DeclareMathOperator{\Aut}{Aut}

\theoremstyle{definition}
\newtheorem{genpos}{Definition}
\newtheorem{Frattini}[genpos]{Definition}

\theoremstyle{remark}
\newtheorem{certify}[genpos]{Remark}

\theoremstyle{plain}
\newtheorem{mmaxdim}[genpos]{Proposition}
\newtheorem{imaxdim}[genpos]{Proposition}
\newtheorem{Frattinidef2}[genpos]{Lemma}
\newtheorem{genFrattini}[genpos]{Lemma}
\newtheorem{mFrattini}[genpos]{Corollary}
\newtheorem{maxFrattini}[genpos]{Lemma}
\newtheorem{Phiofnormalsubgroup}[genpos]{Lemma}
\newtheorem{Phiofpgroup}[genpos]{Lemma}

\newtheorem{counterexample}[genpos]{Proposition}
\newtheorem{goursat}[genpos]{Lemma}
\newtheorem{whiston}[genpos]{Lemma}
\newtheorem{subdirsimple}[genpos]{Lemma}

\newtheorem{S4wrZ3}[genpos]{Proposition}
\newtheorem{SmallGroup}[genpos]{Proposition}
\newtheorem{nilpotent}[genpos]{Proposition}

\theoremstyle{plain}
\newtheorem{ss_def}[genpos]{Definition and Theorem}

\theoremstyle{remark}
\newtheorem{infinite_ss}[genpos]{Remark}

\theoremstyle{plain}
\newtheorem{linearchars}[genpos]{Lemma}
\newtheorem{mofG}[genpos]{Lemma}
\newtheorem{commutemodPhiP}[genpos]{Lemma}
\newtheorem{Kprimeactstrivially}[genpos]{Lemma}
\newtheorem{ss_thm}[genpos]{Theorem}

\newtheorem{GtoK}{Proposition}
\newtheorem{KtoG}[GtoK]{Proposition}

\title{On an Inequality of Dimension-like Invariants for Finite Groups}
\author{Ravi Fernando\thanks{This material is based in part upon work at the 2013 Cornell University Math REU program, supported by the National Science Foundation under Grant No. DMS-1156350.  The author gratefully acknowledges this support, as well as the help and encouragement of Professor R. Keith Dennis.} \\
University of California, Berkeley \\
Berkeley, CA 94720 USA}
\date{\today}
\begin{document}

\maketitle

\begin{abstract}
In this paper, we introduce several notions of ``dimension'' of a finite group, involving sizes of generating sets and certain configurations of maximal subgroups. We focus on the inequality $m(G) \leq \mathrm{MaxDim}(G)$, giving a family of examples where the inequality is strict, and showing that equality holds if $G$ is supersolvable.
\end{abstract}

\section{Introduction and background on generating sets}

For $G$ an arbitrary group, a sequence\footnote{When the ordering does not matter, we will often abuse notation and refer to sequences and sets interchangeably.} $s = (g_1, \dots, g_n)$ of elements of $G$ is said to be a {\bf generating sequence} if we have $\langle g_1, \dots, g_n \rangle = G$.  A sequence $s$, generating or otherwise, is said to be {\bf irredundant} if $\langle g_j : j \neq i \rangle$ is properly contained in $\langle g_i \rangle$ for every $i$.  (The same property is sometimes called {\it independent} or {\it minimal}.)  Then every finite generating sequence of a group contains an irredundant one, since we can simply remove redundant elements one at a time until this is no longer possible.  It is worth noting, however, that strange things can happen in some infinite groups; for example, the reader can check that the additive group $\Q$ has no irredundant generating sets.  However, we will only be interested in the case of finite groups.
\\
\\
Armed with these definitions, we can introduce three notions of ``dimension'' of a finite group $G$ which have been studied extensively.  Let $r(G)$ be the minimum size of a (necessarily irredundant) generating sequence of $G$; let $m(G)$ be the maximum size of an irredundant generating sequence of $G$; and let $i(G)$ be the maximum size of any irredundant sequence in $G$.  (It follows from the definition that $i(G)$ is the maximum of $m(H)$ as $H$ runs over subgroups of $G$.)  Clearly, we have $r(G) \leq m(G) \leq i(G)$.  It is less clear that $m(G) \lneq i(G)$ for some $G$, but examples do exist; we will later give examples of groups satisfying an even stronger inequality than this.
\\
\\
To justify our use of the word ``dimension'', consider the case of the elementary abelian group $G = (\mathbb Z/p \mathbb Z)^n$, which we can view as an $n$-dimensional vector space over $\F_p$.  Here, a generating sequence is just a spanning set, an irredundant sequence is a linearly independent set, and an irredundant generating sequence is a basis.  Since all bases have size $n$, it follows that $r(G) = m(G) = i(G) = n$.  On the other hand, for $H = S_n$ (say, $n > 2$), the reader can find irredundant generating sequences proving that $r(H) = 2$ but $m(H) \geq n-1$.  In fact, it is a nontrivial theorem of Julius Whiston (\cite{whiston})---relying ultimately on the classification of finite simple groups, through the O'Nan-Scott theorem on maximal subgroups of $S_n$---that $m(S_n)$ is exactly $n-1$.  Whiston actually proved much more than this, including that $i(S_n) = n-1$, and that $m(A_n) = i(A_n) = n-2$ for $n > 1$.
\\
\\
Given a finite group $G$, there is an important connection between irredundant generating sequences of $G$ and certain configurations of maximal subgroups $M < G$.  To state this precisely, we first need the following definition.
\begin{genpos}
We say that a family of subgroups $H_i \leq G$, indexed by a set $S$, is in {\bf general position} if it satisfies either of the following equivalent conditions:
\begin{enumerate}
  \item Whenever $\cap_{i \in I} H_i = \cap_{j \in J} H_j$ for $I, J \subseteq S$, we have $I = J$.
  \item For every $i \in S$, the intersection $\cap_{j \neq i} H_j$ properly contains $\cap_{j \in S} H_j$.
\end{enumerate}
To show that (1) implies (2), simply take $I = S$ and $J = S \setminus \{i\}$.  To show the reverse implication (by contrapositive), suppose we are given $I \neq J \subseteq S$ violating (1), and take $i \in I \setminus J$ without loss of generality.  Then we have $\cap_{i \in I} H_i = \cap_{j \in J} H_j$, so
\begin{align}
\cap_{j \in I \cup J} H_j = \cap_{j \in J} H_j = \cap_{j \in J \cup I \setminus \{i\}} H_j.
\end{align}
Intersecting both sides with all $H_k$ for $k \notin I \cup J$ yields $\cap_{j \neq i} H_j = \cap_{\text{all }j} H_j$, contradicting (2).
\end{genpos}
To connect this definition to our main topic, let $s = (g_1, \dots, g_n)$ be an irredundant generating sequence of a group $G$.  Then for each $i$, let $H_i = \langle g_j : j \neq i \rangle$.  We must have $g_i \notin H_i$, since otherwise $H_i$ contains $\langle g_1, \dots, g_n \rangle = G$, contradicting irredundancy.  It follows that the intersection of all $H_i$ contains none of the $g_i$, while the intersection of any $n-1$ of them contains exactly one $g_i$.  In particular, using criterion (2) above, we have shown that the $H_i$ are in general position.
\\
\\
For both theoretical and computational purposes, it is useful to take this argument one step further.  Because $G$ is finite, each of the proper subgroups $H_i < G$ can be enlarged to a maximal subgroup $M_i$.  These $M_i$ contain all $g_j$ with $j \neq i$, but still cannot contain the corresponding $g_i$ by properness, so the same argument shows that they are in general position as well.  Thus, any length-$n$ irredundant generating sequence of a finite group $G$ gives rise to a (possibly nonunique) family of $n$ maximal subgroups of $G$ in general position.  Thus, if we let $\md(G)$ denote the size of the largest family of maximal subgroups of $G$ in general position, we have shown:
\begin{mmaxdim} \label{mmaxdim}
For finite $G$, we have $m(G) \leq \md(G)$.
\end{mmaxdim}
Next, we might ask whether the correspondence can be reversed.  That is, given a family of maximal subgroups of $G$ in general position, can we recover an irredundant generating sequence of the same length?  This is not generally possible; in fact, we will see an example in section 2 where $\md(G)$ is strictly greater than $m(G)$.  However, we can always recover some irredundant (but not necessarily generating) sequence of the same length, which implies:
\begin{imaxdim} \label{imaxdim}
For finite $G$, we have $\md(G) \leq i(G)$.
\begin{proof}
Let $(M_i)_{1 \leq i \leq n}$ be a family of subgroups in general position; they need not even be maximal.  By condition (2) of the definition of general position, we can choose elements $g_i \in (\cap_{j \neq i} M_j) \setminus M_i$ for each $i$.  By construction, we have $g_j \in M_i$ if and only if $j \neq i$.  So for every $i$, the subgroup $\langle g_j : j \neq i \rangle$ is contained in $M_i$ and $\langle \text{all }g_j \rangle$ is not, so the elements $g_1, \dots, g_n$ form an irredundant sequence.  Taking $n = \md(G)$ gives the result.
\end{proof}
\end{imaxdim}
\begin{certify}
Suppose we have a family of subgroups $(H_i)_{i \in S}$ and a family $(g_j)_{j \in S}$ of elements of $G$ indexed by the same set $S$, and suppose that $g_j \in H_i$ holds exactly when $j \neq i$.  Then the argument of Proposition \ref{mmaxdim} shows that the $H_i$ are in general position.  In this case, we say that the $g_j$ {\bf certify} that the $H_i$ are in general position.  We can summarize the last two results as saying that every irredundant generating sequence certifies a family of maximal subgroups in general position, and every such family is certified by some irredundant (but not necessarily generating) sequence.
\end{certify}
Computationally, $\md$ seems to behave more like $m$ than $i$, and it has even been suggested that $\md = m$ in general.  While we will see in the next section that this is false, the connection with maximal subgroups is quite fruitful for computing $m$ for small groups.  Gabriel Frieden has written a program in GAP exploiting this idea.  Roughly speaking, it works by finding all maximal subgroups of a group $G$, looking for large families of them in general position, and then checking whether any of these are certified by an irredundant generating sequence.
\\
\\
Our next definition is particularly important to the theory of generating sets of groups, as we will see immediately and throughout our discussion.
\begin{Frattini}
The {\bf Frattini subgroup} $\Phi(G)$ of a group $G$ is the intersection of all maximal subgroups of $G$.  We say that $G$ is {\bf Frattini-free} if $\Phi(G) = 1$.
\end{Frattini}
\begin{Frattinidef2}
An element $g \in G$ (where $G$ is finite, for convenience but not necessity) lies in $\Phi(G)$ if and only if for every generating set $S$ containing $g$, the set $S \setminus \{g\}$ still generates $G$.  Thus we can say that the Frattini subgroup consists of {\bf non-generators}:  elements that ``contribute nothing to generating $G$''.
\begin{proof}
Let $S$ be any subset of $G$, and let $\langle S \rangle = H \leq G$.  Then $H$ is a proper subgroup of $G$ if and only if it is contained in a maximal subgroup of $G$.  (This is false for infinite $G$:  for example, $\Q$ has no maximal subgroups.)  So $S$ generates $G$ if and only if for every maximal $M < G$, there exists $s \in S$ not in $M$.  It follows that removing an element $g \in \Phi(G)$ does not affect the property of generating $G$.  To prove the converse, suppose $g \in G$ does not belong to some maximal subgroup $M$.  Then the set $S = M \cup \{g\}$ generates a subgroup of $G$ strictly larger than $M$, which must be $G$; but removing $g$ leaves a subset that generates only $M$.  So such a $g$ cannot be a non-generator.
\end{proof}
\end{Frattinidef2}
From the definition, it is clear that $\Phi(G)$ is a subgroup of $G$.  In fact it is a characteristic (and thus normal) subgroup, since any automorphism of $G$ permutes its maximal subgroups and therefore preserves their intersection.  This allows us to take the quotient of $G$ by $\Phi(G)$, which is called the {\bf Frattini quotient}.  Since $\Phi(G)$ is in some sense irrelevant to generating $G$, generation properties of groups are often well-behaved under Frattini quotient, as the next few propositions show.
\begin{genFrattini} \label{genFrattini}
Let $N$ be any normal subgroup of $G$ contained in $\Phi(G)$; for example, $N = \Phi(G)$.  If $\{g_i\}$ is any subset of a finite group $G$, then the $g_i$ generate $G$ if and only if their projections modulo $N$ generate $G/N$.
\begin{proof}
The forward direction is clear.  For the reverse direction, suppose the projections $\overline{g_i}$ generate $G/N$.  Then the larger set $\{g_i\} \cup N$ generates $G$, because every $g \in G$ can be written as the product of a word in the $g_i$ and an element of $N$.  But $N \leq \Phi(G)$ consists of non-generators, so we can remove everything in $N$ from our generating set $\{g_i\} \cup N$ to see that the $g_i$ generate $G$.
\end{proof}
\end{genFrattini}
Notice that we already need $G$ to be a finite group in the lemma above:  if $G = \Q$, then $\Phi(G) = \Q$, because $\Q$ contains no maximal subgroups.  In this case, our argument only allows us to remove finitely many elements of $\Phi(G)$ from a generating set, which is not enough.  Indeed, $\Q/\Phi(\Q) = 0$ is generated by the empty set, and $\Q$ is not even finitely generated. 
\begin{mFrattini} \label{mFrattini}
If $N \trianglelefteq G$ is a normal subgroup contained in $\Phi(G)$, then we have $m(G/N) = m(G)$.  (The same is true of $r(G)$, although we won't use this.)
\begin{proof}
As shown above, quotients by such $N$ do not affect the property of being a generating set.  Since a generating set is irredundant if and only if no proper subset generates $G$, it follows that such quotients also do not affect the property of being an irredundant generating set.  So $G/N$ has an irredundant generating sequence of any given length if and only if $G$ does.
\end{proof}
\end{mFrattini}
Unfortunately, $i(G)$ is not so well-behaved under Frattini quotients.  For example, if $G$ is the wreath product $(\mathbb Z/p\mathbb Z) \wr (\mathbb Z/p\mathbb Z) = (\mathbb Z/p\mathbb Z)^p \rtimes (\mathbb Z/p\mathbb Z)$, where the last factor acts by permuting the previous factors, then it can be shown that $i(G) = m((\mathbb Z/p\mathbb Z)^p) = p$ but $i(G/\Phi(G)) = 2$.  However, it is easy to prove that $\md$ is preserved by Frattini quotients:
\begin{maxFrattini}
If $N$ is a normal subgroup of $G$ contained in $\Phi(G)$, then we have $\md(G/N) = \md(G)$.
\begin{proof}
The subgroups of $G/N$ all have the form $H/N$, where $N \leq H \leq G$, and we have a natural bijection $H \leftrightarrow H/N$ between subgroups of $G$ containing $N$ and subgroups of $G/N$.  It follows that the maximal subgroups of $G$ (which all contain $\Phi(G)$, and thus $N$) correspond to the maximal subgroups of $G/N$.  A family of maximal subgroups $M_i < G$ is in general position if and only if the $M_i/N$ are in general position, so the largest such families have the same size.
\end{proof}
\end{maxFrattini}
We will use two more standard facts about Frattini subgroups; the proofs are left as exercises.
\begin{Phiofnormalsubgroup}
If $N$ is a normal subgroup of $G$, then $\Phi(N) \leq \Phi(G)$.
\end{Phiofnormalsubgroup}
\begin{Phiofpgroup}
If $P$ is a finite $p$-group, then $\Phi(P)$ is the subgroup generated by all commutators and $p$-th powers in $P$.  In particular, $P$ is Frattini-free if and only if it is elementary abelian.
\end{Phiofpgroup}
Notice that using these two lemmas and Corollary \ref{mFrattini}, it is straightforward to calculate $m$ of any finite abelian group $G$.  Specifically, writing $G$ as a direct sum of cyclic groups of prime power order, $G = \oplus_{i = 1}^k \mathbb Z/p_i^{e_i} \mathbb Z$, the reader can show that $m(G) = m(\oplus_{i = 1}^k \mathbb Z/p_i \mathbb Z) = k$ by modding out by each $\Phi(\mathbb Z/p_i^{e_i} \mathbb Z)$.  Notice then that $m(G) \geq m(H)$ for any subgroup $H \leq G$ (which is not true in general), so $i(G) = \max_{H \leq G} (m(H)) = m(G)$ for finite abelian groups $G$.  So finite abelian groups are what we will call \textbf{flat} groups:  groups $G$ satisfying $m(G) = i(G)$.  Notice that the property of flatness is particularly convenient for studying $\md$, since it turns the inequality $m(G) \leq \md(G) \leq i(G)$ into an equality.


\section{A family of groups where $\md \gg m$}

The results of this section were achieved in collaboration with Atticus Christensen.  The family of counterexamples presented here is a descendant of the first known counterexample, in the group $\PSL(3, 2) \wr (\mathbb Z/2\mathbb Z)$, which was discovered by Gabriel Frieden in 2011.
\\
\\
In this section, we will answer the question of whether $\md = m$ holds for all finite groups.  The answer turns out to be emphatically ``no'', in the sense that we can exhibit a family of groups for which $m$ is bounded and $\md$ is not.  Fix $n \geq 5$ and $p$ prime, and let $G$ be the wreath product $A_n \wr \mathbb (Z/p\mathbb Z) = (A_n)^p \rtimes \mathbb (Z/p\mathbb Z)$, where the $\mathbb Z/p\mathbb Z$ in the semidirect product acts by cyclic permutation of the factors.  We claim:
\begin{counterexample}
For $G = A_n \wr (\mathbb Z/p\mathbb Z)$ with $n \geq 5$, we have $\md(G) \geq p(n-3)$ but $m(G) \leq n$.  In particular, for $n = 5$, this gives $\md(G) \geq 2p$ but $m(G) \leq 5$.
\begin{proof}
First we consider $\md(G)$.  Say the $j$th copy of $A_n$ acts on the points $1^{(j)}, \dots, n^{(j)}$, so that the full group $G$ acts on the set $\{i^{(j)} : 1 \leq i \leq n, 1 \leq j \leq p\}$.  One can show (and we will check later) that for any sequence $k_1, \dots, k_p$, the setwise stabilizer of the set $\{k_1^{(1)}, \dots, k_p^{(p)}\}$ is a maximal subgroup of $G$, isomorphic to $A_{n-1} \wr \mathbb Z/p\mathbb Z$.  Now, for $1 \leq j \leq p$ and $2 \leq i \leq n - 2$, let $M_{ij}$ be the setwise stabilizer of $\{1^{(1)}, 1^{(2)}, \dots, i^{(j)}, \dots, 1^{(p)}\}$.  We claim that these $p(n-3)$ maximal subgroups are in general position; in particular, that their intersection is trivial, but the intersection of any $p(n-3)-1$ of them is nontrivial.
\\
\\
Suppose $g$ belongs to the intersection of the $M_{ij}$.  Then $g$ stabilizes the sets $A = \{2^{(1)}, 1^{(2)}, \dots, 1^{(p)}\}$ and $B = \{3^{(1)}, 1^{(2)}, \dots, 1^{(p)}\}$, since $g$ belongs to $M_{21}$ and $M_{31}$ respectively, so $g$ stabilizes $A \setminus B = \{2^{(1)}\}$.  It follows that $g$ has trivial $\mathbb Z/p\mathbb Z$-component, so whenever $g$ stabilizes $\{k_1^{(1)}, \dots, k_p^{(p)}\}$ setwise, it must stabilize each of the individual points.  Since $g$ belongs to each of the $M_{ij}$, we get that $g$ stabilizes $i^{(j)}$ for all $1 \leq i \leq n-2$ and $1 \leq j \leq p$.  The identity permutation is the only such element of $A_n \wr \mathbb Z/p\mathbb Z$, so we have $\cap M_{ij} = 1$ as claimed.  On the other hand, for any choice of $(i_0, j_0)$, the 3-cycle $(i_0^{(j_0)}, (n-1)^{(j_0)}, n^{(j_0)}) \in A_n \wr (\mathbb Z/p\mathbb Z)$ belongs to all $M_{ij}$ except $M_{i_0 j_0}$, so we have shown that the $M_{ij}$ are in general position.  Therefore, we have $\md(G) \geq p(n-3)$ as claimed.
\\
\\
In order to bound $m(G)$, we will need the following two lemmas.  The first is a classical result of Goursat; the second was implicitly used by Whiston in \cite{whiston}, and versions of it were formulated in \cite{cc} and \cite{keen}.
\begin{goursat}
(Goursat)  Suppose $H$ is a subdirect product of two groups $G$ and $G'$; that is, $H$ is a subgroup of $G \times G'$ such that the projections $p_1: H \to G$ and $p_2: H \to G'$ are both surjective.  Let $N'$ and $N$ be the kernels of $p_1$ and $p_2$; these can be identified as normal subgroups of $G'$ and $G$, respectively.  Then there exists an isomorphism $\varphi: G/N \to G'/N'$ such that $H = \{(g, g') \in G \times G': \varphi(\overline{g}) = \overline{g'}\}$.
\end{goursat}
We leave the proof as an exercise.  In fact, we will only be concerned with the case where $G'$ is a simple group.  In this case, $N'$ must equal $G'$ or $1$, so we have a dichotomy between two types of subdirect products.  The case $N' = G'$ gives $N = G$ and $H = G \times G'$.  The case $N' = 1$ gives $\varphi: G/N \stackrel{\sim}{\to} G'$, so $H = \{(g, \phi(g)): g \in G\}$, where $\phi: G \to G'$ is a surjective homomorphism given by $\phi(g) = \varphi(\overline g)$.
\begin{whiston}
(Whiston)  Suppose $(g_1, \dots, g_m)$ is an irredundant generating sequence for some group $G$, and $N \trianglelefteq G$ is a normal subgroup.  Then, possibly after reordering the $g_i$, there exists some $k \leq n$ and some elements $h_{k+1}, \dots, h_m \in N$ such that the projections $\overline{g_1}, \dots, \overline{g_k}$ form an irredundant generating sequence of $G/N$ and $g_1, \dots, g_k, h_{k+1}, \dots, h_m$ form a new irredundant generating sequence for $G$.
\begin{proof}
Since $g_1, \dots, g_m$ generate $G$, their projections generate $G/N$, so we can remove some elements until we have an irredundant generating sequence, which we call $(g_1, \dots, g_k)$ after reordering.  Because the projections $\overline{g_1}, \dots, \overline{g_k}$ generate $G/N$, we can find for each $i > k$ some $x_i \in \langle g_1, \dots, g_k \rangle$ such that $g_i x_i \in N$.  So let $h_i = g_i x_i$.  Then the elements $g_1, \dots, g_k, h_{k+1}, \dots, h_m$ generate $G$, because they generate all of the original $g_i$ via the identities $g_i = h_i x_i^{-1}$ for $i > k$.  So we only need to show that no proper subset of $\{g_1, \dots, g_k, h_{k+1}, \dots, h_m\}$ generates $G$.  To prove this, first note that if we remove any of $g_1, \dots, g_k$, then the projections no longer generate $G/N$, because $g_1, \dots, g_k$ form an irredundant sequence and $h_i \in N$.  Suppose on the other hand that some $h_i$ is not needed in our generating sequence, so we can write $h_i$ as a word in $g_1, \dots, g_k$ and the $h_j$'s with $j \neq i$.  Expanding each $h_j$ as $g_j x_j$ (with $x_j \in \langle g_1, \dots, g_k \rangle$), we get an expression for $g_i = h_i x_i^{-1}$ in terms of $\{g_j: j \neq i\}$.  (In particular, recall that $g_i$ was not used to construct any of the $x_j$, because $i > k$.)  This contradicts the irredundancy of our original sequence, so the new sequence $(g_1, \dots, g_k, h_{k+1}, \dots, h_m)$ must indeed be an irredundant generating sequence of $G$.
\end{proof}
\end{whiston}
Now we are ready to show that $m(G) \leq n$, where $G = A_n \wr (\mathbb Z/p\mathbb Z)$.  Let $g_1, \dots, g_m$ be an irredundant generating sequence for $G$.  Applying Whiston's lemma with $N = A_n^p$, we can assume without loss of generality that $g_1$ generates the quotient $G / A_n^p \cong \mathbb Z/p\mathbb Z$, and all other $g_i$ belong to $A_n^p$.  Next, we claim that after renumbering, we can force $\langle g_1, g_2, \dots, g_{k+1} \rangle \cap A_n^p$ to project onto the entire first coordinate $A_n^p / A_n^{p-1} = A_n$ for some $k \leq n-2$.  If this is so, then conjugating by powers of $g_1$ will force the same subgroup $\langle g_1, g_2, \dots, g_{k+1} \rangle \cap A_n^p$ to project onto the entire $i$th coordinate of $A_n^p$ for all $i$, which will put us not far from generating the whole group.
\\
\\
To prove this claim, begin by taking any family of elements $h_1, \dots, h_k \in A_n^p$.  With $g_1$ as above, the group $\langle g_1, h_1, \dots, h_k \rangle \cap A_n^p$ is exactly the group generated by $g_1^p$ and the conjugates of all $h_i$ by powers of $g_1$.  It follows that $\langle g_1, h_1, \dots, h_k \rangle \cap A_n^p$ surjects onto the first coordinate of $A_n^p$ if and only if $A_n$ is generated by the first coordinate of $g_1^p \in A_n^p$ and the first coordinates of $g_1^t h_i g_1^{-t}$, for all $t$ (in fact, $0 \leq t < p$ suffices) and $1 \leq i \leq k$.  By assumption, this holds for the family $(h_1, \dots, h_k) = (g_2, \dots, g_m)$, since we were given a generating sequence of $G$.  But we know that $m(A_n) = n-2$, so we can choose $k \leq n-2$ elements from the set $\{g_1^p, g_1^t g_i g_1^{-t}\}$ whose first coordinates still generate $A_n$.  These must arise from at most $n-2$ different generators $g_i$ ($i > 1$), so we can indeed choose $g_2, \dots, g_{k+1}$ with $k \leq n-2$ and $\langle g_1, \dots, g_{k+1} \rangle \cap A_n^p$ surjecting onto the first coordinate of $A_n^p$, proving the claim.
\\
\\
We now have an irredundant generating sequence of $G$ whose first $k+1 \leq n-1$ elements generate a subgroup $H$ such that $H$ surjects onto $\mathbb Z/p\mathbb Z$ and $H \cap A_n^p$ surjects onto the first coordinate of $A_n^p$.  By conjugating by powers of $g_1$, then, $H \cap A_n^p$ surjects onto every coordinate of $A_n^p$, so $H \cap A_n^p$ is a subdirect product of $p$ copies of $A_n$.  We now need one more lemma.
\begin{subdirsimple} \label{subdirsimple}
If $S$ is a nonabelian simple group and $K$ is a subdirect product of $r$ copies of $S$ such that for every pair of indices $1 \leq i < j \leq r$, the projection $\pi_{i,j}: K \to S^2$ onto the $i$th and $j$th coordinates is surjective, then $K = S^r$.
\begin{proof}
This is trivially true for $r \leq 2$.  Suppose for the sake of induction that it holds for some $r$, and consider $K \leq S^{r+1}$ satisfying the hypothesis above.  By the inductive hypothesis, $K$ surjects onto the first $r$ coordinates of $S^{r+1}$, so $K$ is a subdirect product of $S^r$ and $S$.  By Goursat's lemma and simplicity of $S$, this must be either $S^{r+1}$ or a fiber product of the form $\{(s_1, \dots, s_r, \varphi(s_1, \dots, s_r))\}$ for some homomorphism $\varphi: S^r \twoheadrightarrow S$.  In the latter case, $\ker \varphi$ is a normal subgroup of $S^r$, and one can easily show that every normal subgroup of $S^r$ is a direct product of a subset of the factors.  Since $S^r / \ker \varphi \cong S$, $\ker \varphi$ must be a direct product of $r-1$ of the factors of $S^r$, so $\varphi$ factors through one of the projections $\pi_i: S^r \to S$.  In other words, we have shown $\varphi(s_1, \dots, s_r) = \alpha(s_i)$ for some $i$ and some $\alpha \in \Aut(S)$.  Then $\pi_{i,r+1}(K)$ consists only of elements of the form $(s, \alpha(s))$, so $\pi_{i,r+1}$ is not surjective.  This contradicts our assumption, so we must have $K = S^{r+1}$ as desired.
\end{proof}
\end{subdirsimple}
By the same Goursat's lemma argument, the only subdirect products of two copies of $S$ are $S^2$ and subgroups of the form $\{(s, \alpha(s))\}$ for $\alpha \in \Aut(S)$.  We say the latter groups are of {\it diagonal type}.  Thus, if $K$ is a subdirect product of $p$ copies of $S$, then we can describe any pair of coordinates as either {\it independent}, meaning that the projection $\pi_{i,j}: K \to S^2$ is surjective, or {\it diagonally linked}, meaning that $\pi_{i,j}(K)$ is of diagonal type.  If we additionally define each coordinate to be diagonally linked to itself, then it is easy to check that diagonal linkedness is an equivalence relation.
\\
\\
Now take $S = A_n$ and $K = H' \cap A_n^p$, where $H' \leq G$ is some subgroup containing $H = \langle g_1, \dots, g_{k+1} \rangle$ as above.  In particular, we showed that $H \cap A_n^p$ is a subdirect product of $p$ copies of $A_n$, so the possibly larger group $K$ is as well.  Since diagonal linkedness is an equivalence relation, it partitions the coordinates of $A_n^p$ into equivalence classes.  Now consider the effect of conjugating $K$ by $g_1$, which permutes the $p$ copies of $A_n$ nontrivially.  (Notice that $g_1$ normalizes $K$, because $K$ is the intersection of the subgroup $H' \ni g_1$ and the normal subgroup $A_n^p$.)  An easy computation shows that conjugation by $g_1$ is given by permuting coordinates according to the $\mathbb Z/p\mathbb Z$-part of $g_1$ and then conjugating by an appropriate element of $A_n^p$.  Since conjugation by elements of $A_n^p$ does not affect diagonal linkedness, it follows that the diagonal linkedness relation is invariant under a cyclic permutation of the coordinates.  Since $p$ is prime, then, it can only be that all coordinates are diagonally linked or all coordinates are independent.  In the latter case, we have $K = A_n^p$ by Lemma \ref{subdirsimple}, so $H' = \langle K, g_1 \rangle = G$.  It follows that in the former case, $H'$ must be a maximal subgroup of $G$, since enlarging it to $H''$ will yield $H'' \cap A_n^p \gneq K$ and thus $H'' \cap A_n^p = A_n^p$.  In particular, the subgroup $H = \langle g_1, \dots, g_{k+1} \rangle$ that we constructed is either all of $G$ or a maximal subgroup, so we cannot add more than one additional generator without losing irredundancy.  It follows that $m(G) \leq k+2 \leq (n-2)+2 = n$, as claimed.
\end{proof}
\end{counterexample}
Finally, for completeness, we give the proof that the subgroups $M_{ij}$ used to estimate $\md(G)$ were indeed maximal.  Notice that all setwise stabilizers of sets of the form $\{k_1^{(1)}, \dots, k_p^{(p)}\}$ are conjugate to each other by elements of $A_n^p$, so it suffices to consider the case where $k_1 = \cdots = k_p = n$, which yields a subgroup naturally isomorphic to $A_{n-1} \wr (\mathbb Z/p\mathbb Z)$.  It is well-known that $A_{n-1} < A_n$ is maximal, being the stabilizer of a point in a doubly transitive group.  So now we claim that the naturally embedded copy of $M \wr (\mathbb Z/p\mathbb Z)$ in $S \wr (\mathbb Z/p\mathbb Z)$ is always maximal for $S$ nonabelian simple and $M < S$ maximal.  To prove this, take any $g \notin M \wr (\mathbb Z/p\mathbb Z)$, and consider $H = \langle M \wr (\mathbb Z/p\mathbb Z), g \rangle$.  By multiplying $g$ by an appropriate element of $\mathbb Z/p\mathbb Z$, we can obtain an element of $H \cap S^p$ that does not belong to $M^p$, which implies that $H \cap S^p$ surjects onto $S$ in some coordinate.  But $H$ contains nontrivial permutations of the coordinates, and conjugating by these gives us that $H \cap S^p$ surjects onto every coordinate, so it is a subdirect product of $p$ copies of $S$.  Now take indices $1 \leq i < j \leq p$, and consider the projection $\pi_{i,j}$ of $H \cap S^p$ onto its $i$th and $j$th coordinates.  We observed after Lemma \ref{subdirsimple} that the image of $\pi_{i,j}$ must be either $S^2$ or a subgroup of diagonal type.  But since $H$ contains $M^p$, the image of $\pi_{i,j}$ must contain $M^2$, so it can only be $S^2$.  (This uses the fact that $S$ is nonabelian, because the trivial subgroup is maximal in a group of prime order.)  Thus, by Lemma \ref{subdirsimple}, we have $H \cap S^p = S^p$, so $H = S \wr (\mathbb Z/p\mathbb Z)$.  So $M \wr (\mathbb Z/p\mathbb Z)$ plus any other element generates $S \wr (\mathbb Z/p\mathbb Z)$, proving that it is a maximal subgroup.  This completes the proof.


\section{Solvable and nilpotent groups}

Now that we know $\md$ can be much larger than $m$ in general, we turn to the question of what assumptions are needed on $G$ to force $\md = m$.  Our first suspicion might be that the strictness of the inequality for $G = A_n \wr (\mathbb Z/p\mathbb Z)$ may be the result of the many copies of the nonabelian simple group $A_n$ in its composition series, and that $\md = m$ may still hold for solvable groups.  However, a slight variant on our original family of counterexamples dashes our hopes:
\begin{S4wrZ3}
For the solvable group $S_4 \wr (\mathbb Z/3\mathbb Z)$, we have $\md \geq 6$ but $m = 5$.
\begin{proof}
Using the same notation as for $A_n \wr (\mathbb Z/p\mathbb Z)$ above, consider the maximal subgroups $M_{ij}$ for $i = 1, 2, 3$ and $j = 2, 3$.  The proofs that these are maximal and in general position are almost identical to the corresponding proofs for $A_n \wr (\mathbb Z/p\mathbb Z)$, with slight changes because we are working with $S_n$ instead of $A_n$.  The computation $m = 5$ was done in GAP.
\end{proof}
\end{S4wrZ3}

Furthermore, extensive human-assisted computations in GAP, performed by the author and R. Keith Dennis, gave the following result:

\begin{SmallGroup}
The solvable group $G$ listed as {\tt SmallGroup(720, 774)} in GAP's SmallGroups library is the unique smallest group with $\md \neq m$.  It has $m(G) = 4, \md(G) = 5$, and $i(G) = 6$.
\end{SmallGroup}

A little discouraged by the failure of equality in solvable groups, we turn to a simpler class of groups, hoping for a positive result.  Nilpotent groups grant our wish.
\begin{nilpotent}
If $G$ is a finite nilpotent group, then $\md(G) = m(G)$.
\begin{proof}
Let $G$ be a finite nilpotent group.  Recall that both $m$ and $\md$ are preserved under Frattini quotients, so there is no loss of generality in assuming that $G$ is Frattini-free to begin with.  Now recall that a finite nilpotent group is the direct product of its Sylow subgroups:  $G = P_1 \times \cdots \times P_n$.  Since $\Phi(P_i)$ is contained in $\Phi(G) = 1$ for all $i$, all the Sylow subgroups $P_i$ are Frattini-free.  But a Frattini-free $p$-group is elementary abelian, so $G$ is an abelian group.  We already know that abelian groups are flat, so we have $m(G) \leq \md(G) \leq i(G) = m(G)$, and thus $\md(G) = m(G)$ holds for all finite nilpotent groups.
\end{proof}
\end{nilpotent}
The idea of the proof above was to reduce to the case of Frattini-free nilpotent groups, and then understand the structure of such a group well enough to force $m(G) = i(G)$, which implies $\md = m$ by the inequality $m \leq \md \leq i$.  In the next section, we will follow essentially the same outline, but we will work harder to extend our result to a larger class of finite groups, which lies between nilpotent groups and solvable groups.


\section{A proof for supersolvable groups}

Before proving that $\md = m$ for finite supersolvable groups, we give several equivalent definitions of supersolvability of a finite group.  Notice that definition (1) is a strengthened version of solvability, and that finite nilpotent groups (which are the direct products of their Sylow subgroups) satisfy definitions (1-3) by standard facts on $p$-groups.  Thus, at least for finite groups, we have the implications nilpotent $\implies$ supersolvable $\implies$ solvable.  We leave as an exercise the verification that supersolvability is a ``reasonable'' property of finite groups, in that it is closed under taking subgroups, quotients, and finite direct products.
\begin{ss_def}
A \textbf{finite supersolvable group} is a finite group $G$ satisfying any of the following equivalent conditions:
\begin{enumerate}
\item There exists a \textbf{supersolvable series} for $G$; that is, a chain of subgroups $1 = G_0 < G_1 < \cdots < G_n = G$ with each $G_i$ normal in the full group $G$, and each quotient $G_{i+1}/G_i$ cyclic.
\item There exists a \textbf{strong supersolvable series} for $G$; that is, a supersolvable series in which each quotient $G_{i+1}/G_i$ is cyclic of prime order.
\item There exists a strong supersolvable series for $G$ in which the orders of the quotients $G_{i+1}/G_i$ are primes arranged in decreasing order.
\item Every maximal subgroup $H < G$ has prime index.
\end{enumerate}
\begin{proof}
Trivially, (3) implies (2) implies (1).  To prove (1) implies (2), suppose $G$ has a supersolvable series $G_0 < \cdots < G_n$, and consider any quotient $G_{i+1}/G_i \cong \mathbb Z/m \mathbb Z$, where $m > 1$.  If $m = a b$ is composite, then let $H$ be the subgroup of $G_{i+1}$ containing $G_i$ that corresponds to $a \mathbb Z/m \mathbb Z$ in the quotient.  Since $G_i$ and $G_{i+1}$ are normal in $G$, any conjugate of $H$ is a subgroup of $G_{i+1}$ containing $G_i$.  But such a subgroup corresponds to a subgroup of $\mathbb Z/m \mathbb Z$, and so is uniquely determined by its size.  It follows that $H$ is normal in $G$.  Then we have $G_{i+1}/H \cong \mathbb Z/a \mathbb Z$ and $H/G_i \cong \mathbb Z/b \mathbb Z$, so we have lengthened the supersolvable series.  This process can be repeated until all quotients have prime order.
\\
\\
The implication $(2) \implies (3)$ is the finite case of Theorem 2.3 in \cite{pinnock}, attributed to Guido Zappa.  The idea of the proof is to ``switch'' two adjacent factors at a time, using the fact that if $|G_i/G_{i-1}| = p < q = |G_{i+1}/G_i|$, then the group $G_{i+1}/G_{i-1}$ of order $pq$ has a characteristic subgroup of order $q$, which equals $G_i'/G_{i-1}$ for some appropriately chosen $G_i' \triangleleft G$.
\\
\\
The equivalence of (1-3) and (4) takes some more work; this is a theorem of Huppert, and a proof can be found in \cite{hall}, pp. 161-3.  See Theorem 4.23 in \cite{conrad} for a few more equivalent conditions.
\end{proof}
\end{ss_def}
\begin{infinite_ss}
For infinite groups, the conditions above are not all equivalent---indeed, conditions (2) and (3) cannot hold as stated---and only condition (1) is taken as the definition of supersolvability.  Moreover, in the infinite case it is no longer even true that abelian groups are supersolvable; for example, one can check that all supersolvable groups are finitely generated, which rules out groups such as $\Q$.
\end{infinite_ss}
The rest of this section will be spent proving that for all finite supersolvable groups $G$, we have $\md(G) = m(G)$.  To show this, we will prove a stronger statement:  if $G$ is supersolvable with $\Phi(G) \cap G' = 1$, then $G$ is flat; that is, $i(G) = m(G)$.  Before beginning the proof of this, let's see why the claim about $\md$ would follow.  If $G$ is a finite supersolvable group, then the Frattini quotient $H = G/\Phi(G)$ is Frattini-free, and in particular satisfies $\Phi(H) \cap H' = 1$.  So we must have $i(H) = m(H)$, and thus $\md(H) = m(H)$ by the inequality $m \leq \md \leq i$.  But both $m$ and $\md$ are invariant under modding out by Frattini subgroups, so it follows that $\md(G) = m(G)$ as well.
\\
\\
Notice that the claim is stronger than what we actually need.  In particular, it would suffice to prove $i = m$ for the smaller class of Frattini-free supersolvable groups.  The strange-looking condition $\Phi(G) \cap G' = 1$ arose while trying to do just this.  Given a Frattini-free supersolvable group $G$ with $i(G) > m(G)$, we originally tried to construct a proper subgroup $K < G$ satisfying the same conditions.  It turned out that the best we could do was to show $\Phi(K) \cap K' = 1$ instead of $\Phi(K) = 1$; the ``Previous progress'' section gives a more precise statement of why this was the ``best possible'' result.  But with a little more work, it is possible to reach the same conclusion using only the weaker assumption $\Phi(G) \cap G' = 1$, and this allows us to complete our proof by infinite descent.
\\
\\
Besides the results we have already collected, we will use three outside facts.  First, to compute $m$ and $i$, we will rely heavily on Proposition 3.5.1 from \cite{collins}, which tells us that if $N \trianglelefteq G$ is minimal normal and abelian, then $m(G) = m(G/N)$ if $N \leq \Phi(G)$, and $m(G) = m(G/N) + 1$ otherwise.  Second, we will repeatedly make use of the existence and conjugacy of Hall subgroups in finite solvable groups; for example, see Exercise 6.1.33 in \cite{df}.  Third, we will use Maschke's theorem from representation theory, specifically for characteristic-$p$ representations of a finite group whose order is not divisible by $p$.
\\
\\
Our proof begins by studying the structure of finite supersolvable groups a little further, then specializing to the Frattini-free case in order to use representation theory.  Let $G$ be a finite supersolvable group.  Let $p$ be the largest prime dividing $|G|$, and say $|G| = m p^r$, where $p \nmid m$.  By definition (3) above, $G$ has a supersolvable series $1 = G_0 < \cdots < G_k = G$, where all $G_i$ are normal in $G$, and the orders of the quotients are primes in decreasing order.  Then $P = G_r$ is a normal Sylow $p$-subgroup of $G$.  (This is the beginning of what is called a \textit{Sylow tower} of $G$.)  Recall that $G$ possesses a Hall $p'$-subgroup; that is, a subgroup $K$ with $|K| = |G|/|P|$.  Such a subgroup is necessarily a complement to $P$, so that we have $G = P \rtimes K$.  We will keep this notation for the rest of the proof:  when we say $G = P \rtimes K$, we implicitly mean that $G$ is a finite supersolvable group, $P$ is its unique Sylow $p$-subgroup (where $p$ is the largest prime dividing $|G|$), and $K$ is some complement of $P$.
\\
\\
Now suppose additionally that $G$ is Frattini-free.  Then since $P$ is normal, we have $\Phi(P) \leq \Phi(G) = 1$, so $P$ is also Frattini-free.  But Frattini-free $p$-groups are elementary abelian, so $P \cong (\mathbb Z/p \mathbb Z)^r$ for some $r$.  Then conjugation by $K$ gives us a representation $\pi: K \to \Aut(P) = \GL_r(\F_p)$, with $\pi(x) = (v \mapsto x v x^{-1})$ for $x \in K, v \in P$.  Since $p$ does not divide the order of $K$, we can apply Maschke's theorem to see that the characteristic-$p$ representations of $K$ are completely reducible.  This brings us to an important structural lemma.
\begin{linearchars} \label{linearchars}
If $G = P \rtimes K$ is a finite supersolvable group with $P \cong (\mathbb Z/p\mathbb Z)^r$ (for example, if $G$ is Frattini-free) and $\pi: K \to \Aut(P)$ is the representation given by conjugation, then $\pi$ decomposes into linear characters.  
\begin{proof}
Returning to the supersolvable series from which we constructed $P$, there exists a series of subgroups $1 = P_0 < P_1 < \cdots < P_r = P$, all normal in the full group $G$, with $|P_{i+1}/P_i| = p$ for each $i$.  Viewing $P$ as a $K$-module via $\pi$, then, the $P_i$ form a chain of $K$-submodules of $P$.  By complete reducibility, we can write $P_{i+1} = P_i \oplus Q_{i+1}$ for some one-dimensional submodules $Q_1, \dots, Q_r \leq P$.  Then $P$ is the direct sum of the $Q_i$, as desired.
\end{proof}
\end{linearchars}
So we have shown that if $G = P \rtimes K$ is Frattini-free and supersolvable, then $P \cong (\mathbb Z/p\mathbb Z)^r$ and $K$ acts separately on its coordinates, for some choice of basis.  We say that such a basis \textit{diagonalizes the action of $K$}.
\\
\\
Next, with $G = P \rtimes K$ still Frattini-free and supersolvable (and $P$ still a normal Sylow $p$-subgroup), we claim that $m(G) = r + m(K)$.  To prove this, first notice that we can obtain $K = G/P$ by beginning with $G$ and repeatedly modding out by (at most $r$) minimal normal subgroups, all of which will be abelian.  Thus, by repeated application of Proposition 3.5.1 in \cite{collins}, we get that $m(G) \leq m(K) + r$.  But we can easily exhibit an irredundant generating set of $G$ of size $r + m(K)$:  take an irredundant generating set of $K$ of size $m(K)$, and append to it any basis of $P$ that diagonalizes the action of $K$.  This generates $G$ because it generates both $K$ and $P$, but throwing out any generators from $K$ will make it impossible to generate the quotient $G/P = K$, and throwing out a basis vector from $P$ will make it impossible to generate any nonzero entry in the corresponding coordinate of $P$.  So we have $r + m(K) \leq m(G) \leq r + m(K)$, giving equality.
\\
\\
Moreover, the formula above can be written as $m(G) = m(P) + m(K)$, which is even true if $G$ isn't Frattini-free.  To prove this, recall that quotienting a group by any normal subgroup contained in its Frattini subgroup does not change the value of $m$, so the identity we just showed gives $m(G) = m(G/\Phi(P)) = m((P/\Phi(P)) \rtimes K) = m(P/\Phi(P)) + m(K) = m(P) + m(K)$.  (This once again uses the fact that $\Phi(P) \leq \Phi(G)$ holds for normal subgroups $P \trianglelefteq G$.  The group $(P/\Phi(P)) \rtimes K$ makes sense because $\Phi(P)$ is a characteristic subgroup of $P$, which must be preserved by $K$.)  So we have shown:
\begin{mofG} \label{mofG}
For a finite supersolvable group $G = P \rtimes K$, we have $m(G) = m(P) + m(K)$.
\end{mofG}
Now suppose $G$ is supersolvable with $\Phi(G) \cap G' = 1$, and $i(G) > m(G)$.  We will eventually apply infinite descent by showing that $K$ satisfies the same hypotheses.  First, we will show that if $K$ is flat---that is, if $i(K) = m(K)$---then $G$ is too.  To do this, notice that $\Phi(P) \cap P'$ is contained in $\Phi(G) \cap G'$, and is therefore trivial.  In fact, the Frattini subgroup of a $p$-group contains the commutator subgroup, so we have $P' = \Phi(P) \cap P' = 1$; that is, the $p$-group $P$ is abelian.  (Compare this to the case $\Phi(G) = 1$, in which we proved earlier that $P$ is \textit{elementary} abelian.)  Now let $H$ be any subgroup of $G$, and let $Q$ and $L$ be a Sylow $p$-subgroup and a Hall $p'$-subgroup of $H$, respectively.  By Sylow's theorem and normality of $P$, we have $Q \leq P$, and by a corresponding theorem for Hall subgroups, $L$ is contained in a conjugate of $K$.  Conjugating $H$ appropriately, we can force $L \leq K$, with $Q$ still contained in $P$.  Then we have $H = Q \rtimes L$.  By Lemma \ref{mofG}, it follows that $m(H) = m(Q) + m(L)$.  But since $Q \leq P$ and $L \leq K$, this is bounded by $i(P) + i(K)$.  Since $P$ is abelian, it is flat; that is, $i(P) = m(P)$.  So if $K$ is also flat, then we have $m(H) = m(Q) + m(L) \leq m(P) + m(K) = m(G)$; that is, $m$ of any subgroup of $G$ is bounded by $m(G)$.  So we have shown that if $\Phi(G) \cap G' = 1$, then $K$ flat implies $G$ flat.  By contrapositive, if $i(G) > m(G)$, then $i(K) > m(K)$ as well.
\\
\\
If there exists a finite supersolvable group $G = P \rtimes K$ with $\Phi(G) \cap G' = 1$ but $i(G) > m(G)$, then the same is true of $K$, except possibly the condition $\Phi(K) \cap K' = 1$.  (Recall in particular that subgroups of supersolvable groups are supersolvable.)  In order to apply infinite descent, we must show that this is actually the case.  We will accomplish this by studying the conjugation action $\pi$ of $K$ on $P$ more closely; the main step will be to prove Lemma \ref{Kprimeactstrivially}, that the commutator subgroup $K'$ acts trivially.  To accomplish this, we will need the following easy lemma:
\begin{commutemodPhiP} \label{commutemodPhiP}
For a finite supersolvable group $G = P \rtimes K$ with $\Phi(G) \cap G' = 1$, if $g \in K$ and $v \in P$ commute modulo the normal subgroup $\Phi(P)$, then $g$ and $v$ commute in $G$.
\begin{proof}
Take $g \in K$ and $v \in P$, and suppose the commutator $[g, v] = g v g^{-1} v^{-1}$ belongs to $\Phi(P)$.  Since $\Phi(P) \leq \Phi(G)$, we have $[g, v] \in \Phi(G)$.  But since this is a commutator, it is also in $G'$.  So we have $[g, v] \in \Phi(G) \cap G' = 1$, proving that $g$ and $v$ commute in $G$.
\end{proof}
\end{commutemodPhiP}
Writing the abelian $p$-group $P$ additively, and using dot notation for the action of $K$, we can express this by saying that $g \cdot v - v \in \Phi(P)$ implies $g \cdot v = v$.  This situation will arise in the next lemma, which we are now ready to tackle.
\begin{Kprimeactstrivially} \label{Kprimeactstrivially}
For a finite supersolvable group $G = P \rtimes K$ with $\Phi(G) \cap G' = 1$, the commutator subgroup $K'$ acts trivially on (that is, centralizes) $P$.
\begin{proof}
Consider the action of $K$ on $P/\Phi(P)$.  Since $P/\Phi(P)$ is elementary abelian, applying Lemma \ref{linearchars} to $G/\Phi(P) = (P/\Phi(P)) \rtimes K$ shows that there exists a basis of $P/\Phi(P)$ diagonalizing this action.  (Notice that although $G/\Phi(P)$ may not be Frattini-free, we can still apply Lemma \ref{linearchars} because its Sylow $p$-subgroup $P/\Phi(P)$ is elementary abelian.)  Let $\overline v_1, \dots, \overline v_r$ be such a basis, and lift the $\overline v_i$ to elements $v_i \in P$; these form a generating set for $P$ by Lemma \ref{genFrattini}.  Now let $v$ be one of the $v_i$, and let $g, h \in K$.  We will compare the actions of $gh$ and $hg$ on $v$.  Note that by construction of $v$, we have $g \cdot \overline v = c \overline v$ (in additive group notation) for some integer $c$, and similarly $h \cdot \overline v = d \overline v$.  Lifting to $P$, we have $g \cdot v = c v + p x$ for some $p x \in \Phi(P)$, and similarly $h \cdot v = d v + p y$.  Expanding using the homomorphism property of the action gives:
\begin{align}
h g \cdot v & = h \cdot (c v + p x) = h \cdot c v + h \cdot p x \\
& = c (d v + p y) + h \cdot p x = c d v + p (c y + h \cdot x),
\end{align}
and
\begin{align}
g h \cdot v & = g \cdot (d v + p y) = g \cdot d v + g \cdot p y \\
& = d (c v + p x) + g \cdot p y = c d v + p (d x + g \cdot y).
\end{align}
(In fact, one can check using Lemma \ref{commutemodPhiP} that $h \cdot px = px$ and $g \cdot py = p y$, but we won't need this.)  Notice that the commutator $[g, h] = gh(hg)^{-1}$ sends $hg \cdot v$ to $gh \cdot v$.  These differ by a multiple of $p$, which lies in $\Phi(P)$, so Lemma \ref{commutemodPhiP} implies that we must have $gh \cdot v = hg \cdot v$.  Since the elements $v = v_i$ generate $P$, it follows that $gh$ and $hg$ act identically on all of $P$, and therefore $[g, h]$ acts trivially.  So we have proved that all commutators $[g, h] \in K'$ act trivially on (that is, centralize) $P$, and therefore the same is true of all of $K'$.
\end{proof}
\end{Kprimeactstrivially}
Now let's put all the pieces together.
\begin{ss_thm}
If $G$ is a finite supersolvable group with $\Phi(G) \cap G' = 1$, then $m(G) = i(G)$.  As a corollary, it follows that $\md = m$ for all finite supersolvable groups.
\begin{proof}
Write $G = P \rtimes K$ as before; we first claim that $\Phi(K) \cap K' \leq \Phi(G) \cap G' = 1$.  (Notice that since $K$ is not necessarily normal in $G$, it may not be true that $\Phi(K)$ is contained in $\Phi(G)$, which would make the statement trivial.)  Suppose $g$ belongs to $\Phi(K) \cap K'$.  In particular, by Lemma \ref{Kprimeactstrivially}, $g$ centralizes $P$.  Let $H$ be any maximal subgroup of $G$, and recall that $H$ is conjugate to $Q \rtimes L$ for some subgroups $Q \leq P, L \leq K$.  In fact, if we write $xHx^{-1} = Q \rtimes L$, we can take $x \in P$ without loss of generality, since conjugating by $K$ does not affect the condition that $Q \leq P$ and $L \leq K$.  Recall from definition (4) of finite supersolvable groups that $H$ must have prime index, so either $L = K$ or $L$ is maximal (with prime index) in $K$.  Since $g$ belongs to $\Phi(K) = \bigcap_{M < K \text{ maximal}} M$, it must belong to $L$ in both cases.  Since $g$ centralizes $P$, $g$ must furthermore belong to all $P$-conjugates of $L$.  So $g$ belongs to all maximal subgroups $H = x^{-1} (Q \rtimes L) x < G$, and thus $g \in \Phi(G)$.  But of course $g \in G'$, because $g \in K' \leq G'$, so in fact $g$ belongs to $\Phi(G) \cap G'$.  Thus we have shown that $\Phi(K) \cap K' \leq \Phi(G) \cap G'$, so the former is trivial if the latter is.
\\
\\
Now we can apply infinite descent.  If $G = P \rtimes K$ is a finite supersolvable group with $\Phi(G) \cap G' = 1$ and $i(G) > m(G)$, then we have shown that the supersolvable subgroup $K$ also satisfies $i(K) > m(K)$ (after Lemma \ref{mofG}) and $\Phi(K) \cap K' = 1$ (just now).  But since $G$ is not the trivial group, its Sylow subgroup $P$ is nontrivial, so $K$ is strictly smaller than $G$.  Thus, by infinite descent on the order of $G$, it follows that there is no such group $G$, completing the proof.
\end{proof}
\end{ss_thm}


\section{Previous progress}

Consider three properties of a finite group $G$:
\begin{enumerate}
\item $G$ is a supersolvable group with $\md > m$.
\item $G$ is a Frattini-free supersolvable group with $i > m$.
\item $G$ is supersolvable with $\Phi(G) \cap G' = 1$, and $i > m$.
\end{enumerate}
We conjecture that no finite group satisfies any of these three properties.  (This is confirmed in the previous section.)  Since $m \leq \md \leq i$, where the former two are invariant under modding out by Frattini subgroups, the existence of a group $G$ satisfying (1) implies that $G/\Phi(G)$ satisfies (2).  Trivially, (2) implies (3) for any group, but we can do better.  The next two propositions establish a many-to-one correspondence between groups satisfying (2) and (smaller) groups satisfying (3).
\begin{GtoK} \label{GtoK}
Suppose $G$ satisfies (2), and let $p$ be the largest prime dividing $|G|$.  By general theory of supersolvable groups (Corollary 3.2a in \cite{pinnock}) we can write $G$ as $P \rtimes K$, where $P$ is the unique Sylow $p$-subgroup and $K$ is a complement.  Then $P$ is elementary abelian and $K$ satisfies (3).
\begin{proof}
Since $P$ is normal in $G$, we have $\Phi(P) \leq \Phi(G) = 1$.  Since the Frattini quotient of any $p$-group is elementary abelian, we have $P = P/\Phi(P) = (\mathbb Z/p\mathbb Z)^r$ for some $r$.  By normality of $P$, $K$ acts on this vector space by conjugation, giving a characteristic-$p$ representation $\pi$ of $K$.  Since $p$ does not divide the order of $K$, we can apply Maschke's theorem to see that the characteristic-$p$ representations of $K$ are completely reducible.
\\
\\
We claim that $\pi$ decomposes into linear characters.  By Theorem 2.3 in \cite{pinnock} (attributed to Guido Zappa), there exists a series of subgroups $1 = P_0 < P_1 < \cdots < P_r = P$, all normal in the full group $G$, with $|P_{i+1}/P_i| = p$ for each $i$.  Viewing $P$ as a $K$-module via $\pi$, then, the $P_i$ form a chain of $K$-submodules of $P$.  By complete reducibility, we can write $P_{i+1} = P_i \oplus Q_{i+1}$ for some one-dimensional submodules $Q_1, \dots, Q_r < P$, and then $P$ is the direct sum of the $Q_i$, as desired.
\\
\\
From now on, we will view $P$ as $(\mathbb Z/p\mathbb Z)^r$, with $K$ acting separately on the coordinates.  Given that $G = P \rtimes K$ satisfies condition (2), we claim that $K$ satisfies (3).  Since subgroups of supersolvable groups are supersolvable, it suffices to show that $i(K) > m(K)$ and $\Phi(K) \cap K' = 1$.  We will proceed in this order.
\\
\\
First, we claim $m(G) = m(K) + r$ and $i(G) = i(K) + r$, so that $G$ is flat if $K$ is.  Let $H$ be an arbitrary subgroup of $G$, possibly $G$ itself.  We would like to compute $m(H)$.  Let $Q$ be a Sylow $p$-subgroup of $H$, and let L be a Hall $p'$-subgroup of $H$; that is, a subgroup of order $|H|/|Q|$.  Then $Q \leq P$ by Sylow's theorem and normality of $P$ in $G$, and $L$ is contained in a conjugate of $K$ by Hall's theorem.  Since we only care about the isomorphism class of $H$, we may conjugate it so that $L$ is actually contained in $K$; we will still have $Q \leq P$.
\\
\\
By the same argument used above, we can show that $L$ acts separately on the coordinates of $Q$, for some appropriate choice of basis of the elementary abelian group $Q$.  If $B$ is such a basis, then combining $B$ with any maximal irredundant generating sequence of $L$ yields an irredundant generating sequence of $H$ of length $m(L) + \text{rank}(Q)$.  On the other hand, we can obtain $L = H/Q$ from $H$ by modding out by abelian minimal normal subgroups at most rank$(Q)$ times, so applying Proposition 3.5.1 in \cite{collins} repeatedly gives $m(H) \leq m(L) + \text{rank}(Q)$.  We have shown inequalities in both directions, so we have $m(H) = m(L) + \text{rank}(Q)$.  We can use this formula in two ways.  First, setting $H = G$, we have $m(G) = m(K) + \text{rank}(P) = m(K) + r$.  Second, taking upper bounds on $m(L)$ and $\text{rank}(Q)$ gives $i(G) = \max_{H \leq G} m(H) \leq i(K) + r$, and this bound is in fact realized by the subgroup $H = PL$ where $L \leq K$ is chosen with $m(L) = i(K)$.  So we have proved that $m(G) = m(K) + r$ and $i(G) = i(K) + r$.  Since we assumed that $G$ is not flat, we have $i(K) - m(K) = i(G) - m(G) > 0$, so $K$ is not flat either.
\\
\\
Finally, we must show that $\Phi(K) \cap K' = 1$.  We will proceed by contradiction, assuming $\Phi(K) \cap K' \neq 1$ and concluding that $\Phi(G) \neq 1$ as well.  Suppose a nonidentity element $x \in K$ is contained in both $\Phi(K)$ and $K'$.  Since $K$ acts separately on the coordinates of $P$, the map $K \to \Aut(P)$ given by conjugation has image contained in the abelian group $(\Aut(\mathbb Z/p\mathbb Z))^r = ((\mathbb Z/p\mathbb Z)^{\times})^r$.  It follows that commutators act trivially, and in particular $x \in K'$ is centralized by $P$.
\\
\\
Now we claim $x$ lies in all maximal subgroups of $G$.  If $H < G$ is maximal, then by a general fact on supersolvable groups we have $[G : H] = q$ prime.  As before, let $L$ be a Hall $p'$-subgroup of $H$, so that $L$ is conjugate to either $K$ or one of its maximal subgroups.  In particular, since $G = PK$, $L$ is a $P$-conjugate of either $K$ or a maximal subgroup of $K$.  But $x$ belongs to all maximal subgroups of $K$ by assumption, so $x$ belongs to a $P$-conjugate of $L$.  Since we have already shown that $P$ centralizes $x$, it follows that $x$ belongs to $L$, and thus $x \in H$.  So $x \neq 1$ belongs to the intersection $\Phi(G)$ of all maximal subgroups of $G$, contradicting the assumption that $\Phi(G) \neq 1$.
\\
\\
So we have shown that $K$ is supersolvable and non-flat with $\Phi(K) \cap K' = 1$; that is, $K$ satisfies (3).
\end{proof}
\end{GtoK}

\begin{KtoG} \label{KtoG}
Suppose $K$ satisfies (3).  Then there exist infinitely many primes $p$ and groups $G = P \rtimes K = (\mathbb Z/p\mathbb Z)^r \rtimes K$ such that $G$ satisfies (2).
\begin{proof}
Let $K$ be any finite group satisfying (3).  Say the abelian group $K/K'$ is isomorphic to $(\mathbb Z/n_1 \mathbb Z) \times \cdots \times (\mathbb Z/n_k \mathbb Z)$, and let $p$ be any prime that is congruent to 1 modulo all of the $n_i$.  (Dirichlet's theorem guarantees the existence of infinitely many such $p$.  We can choose $p$ greater than all primes dividing $|K|$ if we want to imitate the situation of the first proposition, but this isn't necessary.)  Then each $\mathbb Z/n_i \mathbb Z$ embeds into the cyclic group $\F_p^{\times}$, since $n_i$ divides $p-1$ by assumption.  So we can define $k$ linear characters $\chi_1, \dots, \chi_k:  K \to \F_p^{\times}$, where $\chi_i$ first projects $K$ onto $K/K'$, then projects this onto its $i$th coordinate $\mathbb Z/n_i \mathbb Z$, and finally embeds this in $\F_p^{\times}$.  The direct sum of these $k$ characters is a characteristic-$p$ representation $\pi:  K \to (\F_p^{\times})^k \leq \Aut((\mathbb Z/p\mathbb Z)^k)$.  Moreover, we have $\ker \pi = K'$, because $\pi$ factors through an injective map $K/K' \cong (\mathbb Z/n_1 \mathbb Z) \times \cdots \times (\mathbb Z/n_k \mathbb Z) \to (\F_p^{\times})^k$.  Now let $P = (\mathbb Z/p\mathbb Z)^k$, and let $G$ be the semidirect product $P \rtimes K$, where $K$ acts on $P$ by $\pi$.  We claim that $G$ satisfies property (2).  Three things must be checked:  that $G$ is supersolvable, that $\Phi(G) = 1$, and that $i(G) > m(G)$.  Notice that once we show the first two of these, the last will follow from a step we used to prove the previous proposition:  for a Frattini-free supersolvable group $G = P \rtimes K$, we have $m(G) = m(K) + k$ and $i(G) = i(K) + k$, where $k$ corresponds to $r$ above; and $K$ is non-flat by assumption.
\\
\\
Now we will show that $G$ is supersolvable.  Since $K$ is assumed to be supersolvable, we are given a supersolvable series $1 = K_0 \triangleleft K_1 \triangleleft \cdots \triangleleft K_{\ell} = K$.  Here, each $K_i$ is normal in $K$, and the quotients $K_{i+1}/K_i$ are cyclic; for convenience (and without loss of generality), we take the quotients to be cyclic of prime order.  For $0 \leq j \leq k$, let $P_j$ denote the subspace spanned by the first $j$ coordinates of the vector space $P = (\mathbb Z/p\mathbb Z)^k$.  Since $K$ acts separately on the coordinates of $P$, each $P_j$ is $K$-invariant and thus normal in $G$.  So consider the series:
\[
1 = P_0 \triangleleft P_1 \triangleleft \cdots \triangleleft P_k = P = PK_0 \triangleleft PK_1 \triangleleft \cdots \triangleleft PK_{\ell} = G.
\]
We have seen that the $P_j$ are normal in $G$.  The $PK_i$ are subgroups of $G$ because $P \triangleleft G$, and they are normal since they are normalized by both $P$ and $K$.  Each quotient of consecutive terms has prime order, so this is indeed a supersolvable series for $G$.
\\
\\
Finally, we claim that $G$ is Frattini-free.  To prove this, we will use two types of maximal subgroups of $G$ to show that $\Phi(G) \leq \Phi(K) \cap K'$, which is trivial by assumption.  First, if $L$ is any maximal subgroup of $K$, then $PL$ is a maximal subgroup of $G$ by index considerations.  (Recall that in finite supersolvable groups, we have a convenient criterion for maximality:  a subgroup is maximal if and only if it has prime index.)  Second, by the same reasoning, if $Q$ is a maximal subgroup of $P$ that is $K$-invariant, then $QK$ is a maximal subgroup of $G$.
\\
\\
Intersecting all subgroups of the first type gives $\Phi(G) \leq P \Phi(K)$.  For the second type, recall that each coordinate of $P$ is $K$-invariant, so the sum of any $k-1$ of the coordinates is $K$-invariant.  This yields $k$ maximal subgroups $Q_1, \dots, Q_k$ of $P$, all $K$-invariant, with trivial intersection.  It follows that $\Phi(G)$ is contained in the intersection $\bigcap_i Q_i K = K$.  But since $\Phi(G)$ is normal in $G$, all of its conjugates are contained in $K$ as well.  (That is, $\Phi(G)$ is contained in the core of $K$.)  In particular, if $x \in \Phi(G)$, then $vxv^{-1} \in K$ for all $v \in P$.  But $x$ and $vxv^{-1}$ project to the same element of $G/P$ $(\cong K)$, since $v$ is trivial in this quotient.  Two elements of $K$ that are congruent modulo $P$ are equal, so we have $x = vxv^{-1}$, showing that every $x \in \Phi(G)$ centralizes $P$.  Thus we have $\Phi(G) \leq K \cap C_G(P)$.  By the construction of $G$ as $P \rtimes_{\pi} K$, the subgroup of $K$ centralizing $P$ is precisely the kernel of $\pi$, and we constructed $\pi$ so that its kernel is precisely the commutator subgroup $K'$.  Thus, $\Phi(G)$ is contained in $K'$, and from earlier it is contained in $P \Phi(K)$, so it is contained in the intersection $P \Phi(K) \cap K' = \Phi(K) \cap K'$.  Since we assumed that $K$ satisfies (3), this is trivial.  So $G$ is Frattini-free, completing the proof.
\end{proof}
\end{KtoG}

The propositions above are interesting for a few reasons.  First, they show, quite constructively, that there exist groups satisfying (2) if and only if there exist groups satisfying (3).  Moreover, they focus the search for a possible proof that no such groups exist.  A first idea at such a proof might be to induct on the number of primes dividing a group's order, proceeding from a group $G$ supposedly satisfying (2) to its Hall subgroup $K$.  But Proposition \ref{KtoG} shows that we can only hope to prove that $K$ satisfies the weaker condition (3), since any such group will have $G$ satisfying (2) sitting ``above'' it.  This suggests that we should begin with (3), not (2), as our inductive hypothesis.  The previous section uses exactly this strategy to prove that $\md = m$ for finite supersolvable groups.

\end{document}